\documentstyle[12pt]{article}
\begin{document}
\renewcommand{\thefootnote}{\fnsymbol{footnote}}
\newpage
\pagestyle{empty}
\setcounter{page}{0}
\newcommand{\dotimes}{\stackrel{}{\dot{\otimes}}}
\renewcommand{\thesection}{\arabic{section}}
\renewcommand{\theequation}{\thesection.\arabic{equation}}
\newcommand{\sect}[1]{\setcounter{equation}{0}\section{#1}}
\newfont{\twelvemsb}{msbm10 scaled\magstep1}
\newfont{\eightmsb}{msbm8}
\newfont{\sixmsb}{msbm6}
\newfam\msbfam
\textfont\msbfam=\twelvemsb
\scriptfont\msbfam=\eightmsb
\scriptscriptfont\msbfam=\sixmsb
\catcode`\@=11
\def\Bbb{\ifmmode\let\next\Bbb@\else
\def\next{\errmessage{Use \string\Bbb\space only in math mode}}\fi\next}
\def\Bbb@#1{{\Bbb@@{#1}}}
\def\Bbb@@#1{\fam\msbfam#1}
\newfont{\twelvegoth}{eufm10 scaled\magstep1}
\newfont{\tengoth}{eufm10}
\newfont{\eightgoth}{eufm8}
\newfont{\sixgoth}{eufm6}
\newfam\gothfam
\textfont\gothfam=\twelvegoth
\scriptfont\gothfam=\eightgoth
\scriptscriptfont\gothfam=\sixgoth
\def\frak{\frak@}
\def\frak@#1{{\fam\gothfam{{#1}}}}
\def\frak@@#1{\fam\gothfam#1}
\catcode`@=12
%
%
%
\def\CC{{\Bbb C}}
\def\NN{{\Bbb N}}
\def\QQ{{\Bbb Q}}
\def\RR{{\Bbb R}}
\def\ZZ{{\Bbb Z}}
\def\cA{{\cal A}}          \def\cB{{\cal B}}          \def\cC{{\cal C}}
\def\cD{{\cal D}}          \def\cE{{\cal E}}          \def\cF{{\cal F}}
\def\cG{{\cal G}}          \def\cH{{\cal H}}          \def\cI{{\cal I}}
\def\cJ{{\cal J}}          \def\cK{{\cal K}}          \def\cL{{\cal L}} 
\def\cM{{\cal M}}          \def\cN{{\cal N}}          \def\cO{{\cal O}}
\def\cP{{\cal P}}          \def\cQ{{\cal Q}}          \def\cR{{\cal R}} 
\def\cS{{\cal S}}          \def\cT{{\cal T}}          \def\cU{{\cal U}}
\def\cV{{\cal V}}          \def\cW{{\cal W}}          \def\cX{{\cal X}}
\def\cY{{\cal Y}}          \def\cZ{{\cal Z}}
\def\qed{\hfill \rule{5pt}{5pt}}
\def\id{\mbox{id}}
\def\ggo{{\frak g}_{\bar 0}}
\def\uqggo{\cU_q({\frak g}_{\bar 0})}
\def\uqggp{\cU_q({\frak g}_+)}
\def\half{\frac{1}{2}}
\def\btf{\bigtriangleup}
\newtheorem{lemma}{Lemma}
\newtheorem{prop}{Proposition}
\newtheorem{theo}{Theorem}
\newtheorem{Defi}{Definition}

\vfill
\vfill
\begin{center}

{\LARGE {\bf {\sf 
Jordanian quantum spheres
}}} \\[0.8cm]

{\large   
R. Chakrabarti$^{1}$ and J. Segar$^{2}$}
\begin{center}
{\em 
$^{1}$Department of Theoretical Physics, University of Madras, \\
Guindy Campus, 
Madras-600 025, India. \\
$^{2}$Department of Physics, R. K. M. Vivekananda College,
Madras-600 004, India. \\[2.8cm]
}
\end{center}

\end{center}

\smallskip

\smallskip 

\smallskip

\smallskip

\smallskip

\smallskip 

\begin{abstract}
We introduce and investigate a one parameter family of quantum
spaces invariant under the left (right) coactions of the group-like
element ${\cal T}_{h}^{(j=1)}$ of the Jordanian function algebra 
$Fun_{h}(SL(2))$. These spaces may be regarded as Jordanian 
quantization of the two-dimensional spheres.   

\end{abstract}

\vfill
\newpage 

\pagestyle{plain}

\sect{Introduction}
Noncommutative quantum spheres were first studied by Podles 
~\cite{P87} by considering a one-parameter family of quantum
homogeneous spaces for the standard $q$-deformed quantum
group $SL_{q}(2)$. Exterior algebras of differential forms
on these $q$-deformed quantum spaces were investigated
~\cite{P89},\cite{P92} and classification of covariant
differential calculi on Podles quantum $q$-spheres was
obtained ~\cite{AS94}. Orthogonal bases for the algebra of 
functions on Podles quantum $q$-spheres were determined 
~\cite{NM90} in terms of the big $q$-Jacobi polynomials.
An elegant description of the $q$-spherical functions was
obtained ~\cite{K93} by considering left and right invariances
in the infinitesimal sense with respect to twisted primitive
elements of the $q$-deformed universal enveloping algebra
$U_{q}(sl(2))$. 

On the other hand it is well-known that the Lie group $SL(2)$
admits ~\cite{Ku92} two distinct quantum group deformations 
with central quantum determinant: the standard 
$q$-deformation $Fun_{q}(SL(2))$
and the Jordanian deformation $Fun_{h}(SL(2))$ ~\cite{DMMZ90}.
On the quantum algebra level, the Jordanian deformation 
$U_{h}(sl(2))$ of the corresponding classical enveloping
algebra was constructed~\cite{O92}. The differential
calculas on the Jordanian $h$-deformed quantum plane was
investigated ~\cite{A93}, and the classification of the
bicovariant diffential calculi on the Jordanian function
algebra $Fun_{h}(SL(2))$ was obtained~\cite{JC98}.
Using a contraction technique developed earlier 
~\cite{AKS95},\cite{ACC98}, a general recipe for obtaining
the finite dimensional group-like 
elements ${\cal T}_{h}^{(j)}$ of the function algebra 
$Fun_{h}(SL(2))$ was described in ~\cite{CQ99}.

In the present work our starting point is the explicit construction
of the matrix  ${\cal T}_{h}^{(j=1)}$ obtained in 
~\cite{CQ99}. Investigating the left (right) coaction of 
${\cal T}_{h}^{(j=1)}$ on quantum spaces with three noncommuting
generators we obtain a one parameter family of algebras invariant 
under the coaction. These invariant quantum spaces may be regarded
as Jordanian deformation of the two-dimensional spheres. We 
demonostrate that at the algebraic level there exist an invertible map 
between the quantum spaces invariant under the left and the right 
coactions of the function algebra $Fun_{h}(SL(2))$ respectively. 
Moreover, the quantum spaces may be realized as the subalgebras of the 
function algebra $Fun_{h}(SL(2))$. Lastly, following ~\cite{K93},
the above embedding may be obtained by using the infinitesimal
invariance of the generating elements with respect to a corresponding  
twisted primitive element of the dual Jordanian enveloping 
algebra $U_{h}(sl(2))$.

A motivation for the present work is the recent interest in
formulating field theory on noncommutative spaces, both from the
point of view of finding new insights into the renormalization
problem ~\cite{F96}-\cite{GMS00}; and the realization that such 
noncommutative spaces are induced by certain sectors of string
theory, particularly open strings ending on $D$-branes with a
background $B$ field ~\cite{DH98}-\cite{SW99}. We, therefore, expect
exploration of new noncommutative spaces to be fruitful. 
 
\sect{The left and the right coactions of the ${\cal T}_{h}^{(j=1)}$
matrix}
The defining relations for the Jordanian deformed $U_{h}(sl(2))$
algebra read ~\cite{O92}
\begin{eqnarray}
&&[H , T^{\pm 1 }] = T^{\pm 2} - 1,\nonumber\\ 
&&[H, Y] = -\half ( Y(T+T^{-1}) + (T+T^{-1} )Y ),\nonumber\\ 
&&[T^{\pm 1} , Y ]  =  \pm \frac{h}{2} (H T^{\pm 1} + T^{\pm 1} H).
\label{eq:jalg}
\end{eqnarray}
The corresponding coproduct $(\bigtriangleup)$, counit 
$(\varepsilon)$  and the antipode $(S)$ maps are given by
\begin{eqnarray}
&&\bigtriangleup (T^{\pm 1}) = T^{\pm 1} \otimes  T^{\pm 1}, \,\,
\bigtriangleup(Y) = Y \otimes T + T^{-1} \otimes Y, \,\,
\bigtriangleup(H) = H \otimes T + T^{-1} \otimes H,\nonumber\\
&&\varepsilon (T^{\pm 1} ) = 1,\qquad \varepsilon (Y) = 
\varepsilon (H) = 0,\nonumber\\
&& S (T^{\pm 1}) = T^{\mp 1}, \quad 
S (Y) = -TYT^{-1},\quad S(H) = -THT^{-1}. 
\label{eq:jcoalg}
\end{eqnarray}
The dual Jordanian function algebra $Fun_{h}(SL(2))$ is generated
by the elements of the matrix ${\cal T}_{h}^{(j=\half)}=\left(
\begin{array}{cc}a & b\\c &d\end{array} \right)$, obeying the 
relations
\begin{eqnarray}
&&[a, b] = h \left(a^2 - 1\right), \qquad [a, c] =  - h c^2,
\nonumber\\
&&[b, d] = - h \left(d^2 - 1\right), \qquad [c, d]  =  h c^2,
\nonumber\\
&& [a, d] =  h ac - h dc,\qquad  [b, c] = - h ac - h cd,
\label{eq:Funhalg}
\end{eqnarray}
where the central determinant $D = ad - bc - h ac$ satisfies the
condition
\begin{equation}
D = 1.
\label{eq:uniDet}
\end{equation} 
The coalgebraic properties for the generating elements $(a, b, c, d)$
may be summarized as
\begin{eqnarray}
&&\bigtriangleup\left(\begin{array}{cc}a & b\\c & d\end{array}\right)
= \left(\begin{array}{cc}a & b\\c & d\end{array}\right)
\dotimes \left(\begin{array}{cc}a & b\\c &d\end{array} \right),
\nonumber\\
&&\varepsilon \left(\begin{array}{cc}a & b\\c & d\end{array} \right) 
= \left(\begin{array}{cc}1 & 0\\0 & 1\end{array}\right),\nonumber\\
&&S \left(\begin{array}{cc}a & b\\c & d\end{array}\right)
= \left(\begin{array}{cc}d - h\,c & - b + h\,(a - d) + h^2\,c\\
- c & a + h\,c\end{array}\right).
\label{eq:Funhcoalg}
\end{eqnarray}
A recipe based on a contraction procedure has been provided in
~\cite{CQ99} for obtaining arbitrary finite dimensional group-like
elements ${\cal T}_{h}^{(j)}$ of the $Fun_{h}\left(SL(2)\right)$
algebra. These ${\cal T}_{h}^{(j)}$ matrices may be regarded as the
finite dimensional representations of the dual form between the
Jordanian Hopf algebras $Fun_{h}\left(SL(2)\right)$ and 
$U_{h}(sl(2))$. Using the determinantal condition ~(\ref{eq:uniDet}),
here we reproduce the explicit determination of the
${\cal T}_{h}^{(j=1)}$ obtained in ~\cite{CQ99}: 
\begin{equation}
{\cal T}_{h}^{(j=1)} = \left(
\begin{array}{ccc}
a^2 + { \frac{1}{4}} h^2 c^2 &
2 ab - h a^2 + h + {\half} h^2 cd  &
2 b^2 - {\frac{3}{2}} h^2 a^2 + { \half} h^2 d^2 \\
{} & - { \frac{1}{4}} h^3c^2 & 
+ h^2 - { \frac{3}{8}}h^4c^2\\
{} & {} & {}\\
ac + { \half} h c^2 &
1 + 2 bc + h ac + h cd &
2 bd + h d^2 - h - { \frac{3}{2}} h^2 ac \\ 
{}& - { \half} h^2 c^2 &  -{ \frac{3}{4}} h^3 c^2 \\
{} & {} & {}\\
{\half} c^2 &
cd - {\half} h c^2 &
d^2 - { \frac{3}{4}} h^2 c^2 \\
\end{array}\right).
\label{eq:T-one-h}
\end{equation}
The group-like property of the Jordanian monodromy matrix 
${\cal T}_{h}^{(j=1)}$ is manifest in its coproduct relation
\begin{equation}
\bigtriangleup \left({\cal T}_{h}^{(j=1)}\right) =
{\cal T}_{h}^{(j=1)} \dotimes {\cal T}_{h}^{(j=1)}.
\label{eq:gruprop}
\end{equation}

We define a left coaction ${\sf{\phi}_{L}}$ acting on the 
generators of a quantum space $\left({\sf x_{i}\, | i= -1, 0, 1} 
\right)$ by the relation 
\begin{equation}
{\sf{\phi}_{L}(x_{i}) = \sum _{j=-1}^{1}} 
\left({\cal T}_{h}^{(j=1)} \right)_{\sf i\,j} \otimes {\sf x_{j}} 
\quad \hbox{for} \quad {\sf i = (-1, 0, 1)}.
\label{eq:lcoac}
\end{equation}    
As a consequence of the coproduct property ~(\ref{eq:gruprop})
the coaction ${\sf{\phi}_{L}}$ satisfies the identity
\begin{equation}
\left( \hbox{id} \otimes {\sf{\phi}_{L}} \right) \circ
{\sf{\phi}_{L}(x_{i}}) = \left( \bigtriangleup \otimes \hbox{id}
\right) \circ {\sf{\phi}_{L}(x_{i}})
\quad \hbox{for} \quad {\sf i = (-1, 0, 1)}.
\label{eq:phiLid}
\end{equation} 
For two complex numbers $\left( k, \beta \right)$, where $\beta 
\neq 0$, we denote by ${\cal X}_{(h;\; k, \beta)}$ the
algebra with the noncommuting generators $\left({\sf x_{i}\, 
| i= -1, 0, 1} \right)$ obeying the defining relations: 
\begin{eqnarray}
{\sf x_{0}}^2 - {\sf x_{-1} x_{1} - x_{1} x_{-1}} + 8\, h^2\,
{\sf x_{1}}^2 + 8\, kh^2 {\sf x_{1}} &=&  \beta, \nonumber\\
{\sf x_{-1} x_{0} -  x_{0} x_{-1}} + h\,{\sf x_{-1} x_{1}}
+ 3\,h\,{\sf x_{1} x_{-1}} + \half \, h^2\, 
{\sf x_{0} x_{1}} - \frac{9}{2}\, h^2\, 
{\sf x_{1} x_{0}} & &\nonumber\\ 
- 6\,h^3\, {\sf x_{1}}^2 + 4\,kh\,{\sf x_{-1}} 
- 4\,kh^2\,{\sf x_{0}} - 6\, kh^3\,{\sf x_{1}} &=& 0,\nonumber\\
{\sf x_{-1} x_{1} - x_{1} x_{-1}} + 2\,h\,{\sf x_{0} x_{1}}
+ 2\,h\,{\sf x_{1} x_{0}} + 4\,kh\,{\sf x_{0}} &=&  0, \nonumber\\
{\sf x_{0} x_{1} - x_{1} x_{0}} + 4\,h\,{\sf x_{1}}^2
+ 4\,kh\, {\sf x_{1}} &=& 0.
\label{eq:lftsphr}
\end{eqnarray}
The coaction ${\sf \phi_{L}}$ defined by ~(\ref {eq:lcoac}) preserves 
the relations ~(\ref {eq:lftsphr}) and is an algebra homomorphism.
The homomorphism ${\sf \phi_{L}}: {\cal X}_{(h;\; k, \beta)}\rightarrow
Fun_{h}\left(SL(2)\right) \otimes {\cal X}_{(h;\; k, \beta)}$ endows
the algebra ${\cal X}_{(h;\; k, \beta)}$ with the structure of 
a left $Fun_{h}\left(SL(2)\right)$ comodule. 

In a similar way a right coaction ${\sf \phi_{R}}$ acting on the
generating elements $({\sf y_{i}\,|\,i=\, -1, 0, 1})$ of a right
quantum space for the funtion algebra  $Fun_{h}\left(SL(2)\right)$
may be defined as  
\begin{equation}
{\sf{\phi}_{R}(y_{i}) = \sum _{j=-1}^{1}} 
{\sf y_{j}} \otimes \left({\cal T}_{h}^{(j=1)} \right)_{\sf j\,i} 
\quad \hbox{for} \quad {\sf i = (-1, 0, 1)}\,.
\label{eq:rcoac} 
\end{equation}
The group-like coproduct property~(\ref{eq:gruprop}) now requires
the right coaction ${\sf \phi_{R}}$ to satisfy the identity:
\begin{equation}
\left( {\sf{\phi}_{R}} \otimes \hbox {id} \right) \circ
{\sf{\phi}_{R}(y_{i}}) = \left( \hbox {id} \otimes \bigtriangleup
\right) \circ {\sf{\phi}_{R}(y_{i}})
\quad \hbox{for} \quad {\sf i = (-1, 0, 1)}.
\label{eq:phiRid}
\end{equation}
The noncommuting generating elements $({\sf y_{i}\,|\,i=\, -1, 0, 1})$
of the right quantum space furnish the  algebra
${\cal Y}_{(h;\,k^{\prime}, \beta^{\prime})}$ defined by the relations
\begin{eqnarray}
{\sf y_{0}}^2 - {\sf y_{-1} y_{1} - y_{1} y_{-1}} + 4\, h^2\,
{\sf y_{-1}}^2 - 8\, k^{\prime}\,h^2 {\sf y_{-1}} 
&=& \beta^{\prime}\,, \nonumber\\
{\sf y_{-1} y_{0} - y_{0} y_{-1}} - 4\,h\,{\sf y_{-1}}^2
+ 4\,k^{\prime}\,h\, {\sf y_{-1}} &=& 0,\nonumber\\
{\sf y_{-1} y_{1} - y_{1} y_{-1}} - 2\,h\,{\sf y_{-1} y_{0}}
- 2\,h\,{\sf y_{0} y_{-1}} 
+ 4\,k^{\prime}\,h\,{\sf y_{0}} &=&  0, \nonumber\\
{\sf y_{0} y_{1} -  y_{1} y_{0}} - h\,{\sf y_{-1} y_{1}}
- 3\, h\,{\sf y_{1} y_{-1}} - \frac{3}{2} \, h^2\, 
{\sf y_{-1} y_{0}} - \frac{5}{2}\, h^2\, 
{\sf y_{0} y_{-1}} & &\nonumber\\ 
- 2\,h^3\, {\sf y_{-1}}^2 + 4\,k^{\prime}\,h\,{\sf y_{1}} 
+ 4\,k^{\prime}\,h^2\,{\sf y_{0}} 
+ 2\, k^{\prime}\,h^3\,{\sf y_{-1}} &=& 0.
\label{eq:rhtsphr}
\end{eqnarray}
The parameters $(k^{\prime}, \beta^{\prime})$, where
$\beta^{\prime} \neq 0$, may be in general distinct from the 
parameters $(k, \beta)$ introduced in ~(\ref{eq:lftsphr}). 
The defining relations~(\ref{eq:rhtsphr}) are preserved under the 
coaction ~(\ref{eq:rcoac}), which acts as an algebra homomorphism
mapping ${\cal Y}_{(h;\,k^{\prime}, \beta^{\prime})}\rightarrow
{\cal Y}_{(h;\,k^{\prime}, \beta^{\prime})} \otimes Fun_{h}(SL(2))$.
The left and the right quantum spaces 
${\cal X}_{(h;\,k,\beta)}$ and 
${\cal Y}_{(h;\,k^{\prime}, \beta^{\prime})}$ defined by the relations
~(\ref {eq:lftsphr}) and ~(\ref {eq:rhtsphr}) respectively, may be
mapped to each other at the algebraic level. We will discuss these 
mappings at the end of the present section.

Under the scaling of the generators ${\sf x_{i}} \rightarrow 
{\sf c}^{-1} {\sf x_{i}}$, it follows from ~(\ref{eq:lftsphr})
that the defining relations of ${\cal X}_{(h;\; k, \beta)}$
go into that of ${\cal X}_{(h;\;{\sf c} k,\, {{\sf c}^2}\,\beta)}$.
Consequently these quantum spaces are isomorphic to each other
for any complex number ${\sf c} \neq 0$. All of these equivalent
quantum spaces may be referred to as left spheres. Using this 
isomorphism we may remove one superfluous parameter of the quantum 
space ${\cal X}_{(h;\; k, \beta)}$. Similar argument holds for the 
right quantum space ${\cal Y}_{(h;\,k^{\prime}, \beta^{\prime})}$ 
defined by the relations ~(\ref{eq:rhtsphr}). The scaling freedom may
again be utilized to eliminate one redundant parameter in
${\cal Y}_{(h;\,k^{\prime}, \beta^{\prime})}$.

The algebras ${\cal X}_{(h;\; k, \beta)}$ and ${\cal Y}_{(h;\,
k^{\prime}, \beta^{\prime})}$ may be realized  by embedding them in 
the function algebra $Fun_{h}\left(SL(2)\right)$. The generating 
elements of the  ${\cal X}_{(h;\; k, \beta)}$ algebra may be 
represented using the linear combinations of the matrix elements
of the monodromy matrix ${\cal T}_{h}^{(j=1)}$ given in
~(\ref{eq:T-one-h}):
\begin{eqnarray}
{\sf x_{-1}} &=&  k\,\left( a^2\,+\, \frac{1}{4}\, h^2\, c^2\right)
\,+\,\rho\,\left(2\,ab\,-\,h\,a^2\,+\,h\,+\,\half h^2 cd\,-\,
\frac{1}{4}\,h^3\,c^2\right)\nonumber\\
& &-\,k\,\left( 2\,b^2\,-\,\frac{3}{2}\,h^2\, a^2\,
+\,\half\,h^2\, d^2\,+\,h^2\,
-\,\frac{3}{8}\,h^4\,c^2\right),\nonumber\\
{\sf x_{0}} &=& k\,\left( ac\,+\,\half\,h\,c^2\right)\,
+\,\rho\,\left( 1\,+\,2\,bc\,+\,h\,ac\,+\,h\,cd\,
-\,\half\,h^2\,c^2\right)\nonumber\\
& &-\,k\,\left( 2\,bd\,+\,h\,d^2\,-\,h\,
-\,\frac{3}{2}\,h^2\,ac\,-\,\frac{3}{4}\,h^3\,c^2\right),\nonumber\\   
{\sf x_{1}} &=& \half\,k\,c^2\,+\,\rho\,\left( cd\,
-\,\half\,h\,c^2\right)\,-\,k\,\left(d^2\,
-\,\frac{3}{4}\,h^2\,c^2\right). 
\label{eq:xembed}
\end{eqnarray}
The representation ~(\ref{eq:xembed}) satisfies the defining relations
~(\ref{eq:lftsphr}). The paramater $\beta$ is now given by 
\begin{equation}
\beta = \rho^2 + 2\, k^2.
\label{eq:xrpbt}
\end{equation}  
The representation ~(\ref{eq:xembed}) is particularly simple if
we select previously mentioned scaling parameter 
${\sf c} = \rho^{-1}$ and consider the limiting value 
$\rho \rightarrow \infty$ for a finite value of $k$.
The consequent algebra ${\cal X}_{(h;\; 0, 1)}$ has an 
embedding in the function algebra  $Fun_{h}\left(SL(2)\right)$
{\it \`{a} la} ~(\ref{eq:xembed}) as follows:
\begin{eqnarray}
{\sf x_{-1}} &=&  2\,ab\,-\,h\,a^2\,+\,h\,+\,\half h^2 cd\,-\,
\frac{1}{4}\,h^3\,c^2\nonumber\\
{\sf x_{0}} &=&  1\,+\,2\,bc\,+\,h\,ac\,+\,h\,cd\,
-\,\half\,h^2\,c^2\nonumber\\
{\sf x_{1}} &=&  cd\,-\,\half\,h\,c^2.
\label{eq:xinfembed}
\end{eqnarray}
The group-like coproduct property (\ref{eq:gruprop}) may be used
to prove the following coproduct structure for the embedding 
(\ref{eq:xembed})
\begin{equation}
\bigtriangleup \left( {\cal X}_{(h;\; k, \beta)} \right)
\subseteq Fun_{h}\left(SL(2)\right) \otimes 
{\cal X}_{(h;\; k, \beta)}.
\label{eq:Xcopro}
\end{equation}

Similarly the generators $({\sf y_{i}\,|\,i=\, -1, 0, 1})$ of the
right quantum space ${\cal Y}_{(h;\,k^{\prime}, \beta^{\prime})}$
may also be represented as linear combinations of the matrix 
elements of the group-like construct ${\cal T}_{h}^{(j=1)}$ given 
in~(\ref{eq:T-one-h}). The embedding is implemented as follows:   
\begin{eqnarray}
{\sf y_{-1}} &=& k^{\prime}\,\left(a^2\,+\,\frac{1}{4}
\,h^2\,c^2\right) +\,\rho^{\prime}\,\left( ac\,+\,\half\,h\,c^2\right)
\,-\half\,k^{\prime}\,c^2\, ,\nonumber\\ 
{\sf y_{0}} &=& k^{\prime}\,\left( 2\,ab\,-\,h\,a^2\,+\,h\,
+\,\half\,h^2\,cd\,-\,\frac{1}{4}\,h^3\,c^2\right)\nonumber\\
& &+\,\rho^{\prime}\,\left( 1\,+\,2\,bc\,+\,h\,ac\,+\,h\,cd\,
-\,\half\,h^2\,c^2\right) 
-k^{\prime}\,\left( cd\,-\,\half\,h\,c^2\right) ,\nonumber\\
{\sf y_{1}} &=& k^{\prime}\,\left( 2\,b^2\,-\,\frac{3}{2}\,h^2\, a^2\,
+\,\half\,h^2\, d^2\,+\,h^2\,-\,\frac{3}{8}\,h^4\,c^2\right)\nonumber\\ 
& &+\,\rho^{\prime}\,\left(2\,bd\,+\,h\,d^2\,-\,h\,
-\,\frac{3}{2}\,h^2\, ac\,-\,\frac{3}{4}\,h^3\,c^2\right) 
- \,k^{\prime}\,\left( d^2\,-\, \frac{3}{4}\, h^2\, c^2\right).
\label{eq:yembed}
\end{eqnarray}
The representation ~(\ref{eq:yembed}) obeys the defining algebraic 
properties ~(\ref{eq:rhtsphr}), and yields the following value
of the parameter $\beta^{\prime}$:
\begin{equation}
\beta^{\prime} = {\rho^{\prime}}^{\;2} 
+ 2 (1 -2\,h^2)\,{k^\prime}^{\;2}.
\label{eq:yrpbtpr}
\end{equation}
The scaling argument made earlier in the context of the left sphere
${\cal X}_{(h;\; k, \beta)}$ may also be used here exactly
similarly to investigate the $\rho^{\prime} \rightarrow 
\infty$ limit for a finite value of $k^{\prime}$. The representation
of the relevant space ${\cal Y}_{(h;\; 0, 1)}$ in this limit reads:
\begin{eqnarray}
{\sf y_{-1}} &=& ac\,+\,\half\,h\,c^2,\nonumber\\ 
{\sf y_{0}} &=&  1\,+\,2\,bc\,+\,h\,ac\,+\,h\,cd\,
-\,\half\,h^2\,c^2,\nonumber\\ 
{\sf y_{1}} &=& 2\,bd\,+\,h\,d^2\,-\,h\,
-\,\frac{3}{2}\,h^2\, ac\,-\,\frac{3}{4}\,h^3\,c^2. 
\label{eq:yinfembed}
\end{eqnarray}
For the embedding (\ref{eq:yembed}) we can also obtain  from
(\ref{eq:gruprop}) the coproduct structure of the 
${\cal Y}_{(h;\,k^{\prime}, \beta^{\prime})}$ 
\begin{equation}
\bigtriangleup \left({\cal Y}_{(h;\,k^{\prime}, \beta^{\prime})} 
\right) \subseteq {\cal Y}_{(h;\,k^{\prime}, \beta^{\prime})}
\otimes Fun_{h}\left(SL(2)\right).
\label{eq:Ycopro}
\end{equation}

The quantum spaces ${\cal X}_{(h;\; k, \beta)}$ and 
${\cal Y}_{(h;\,k^{\prime}, \beta^{\prime})}$ are isomorphic to
each other in the sense that an invertible map exists so that 
the defining properties ~(\ref{eq:lftsphr}) and~(\ref{eq:rhtsphr}) 
of the algebras ${\cal X}_{(h;\; k, \beta)}$ and 
${\cal Y}_{(h;\,k^{\prime}, \beta^{\prime})}$ respectively 
are transformed to each other. The map ${\sf{\pi}}$ reads
\begin{eqnarray}
&&{\sf{\pi}} ({\sf x_{-1}}) = {\sf y_{1}} + 2h\, {\sf y_{0}}
+ 4h^2\,{\sf y_{-1}},\quad {\sf{\pi}}({\sf x_{0}}) = {\sf y_{0}} 
+ 2h\,{\sf y_{-1}},\quad {\sf{\pi}}({\sf x_{1}}) 
= {\sf y_{-1}},\nonumber\\ 
&&\qquad\qquad\qquad\qquad\qquad{\sf{\pi}}(k) = - k^{\prime},
\quad {\sf{\pi}}(\beta) = \beta^{\prime}.
\label{eq:xymap}
\end{eqnarray}
The inverse of the map~(\ref{eq:xymap}) may be readily expressed.

\sect{An infinitesimal characterization of the Jordanian 
quantum sphere}
In this section we give a description of the embedding of the 
algebras ${\cal X}_{(h;\; k, \beta)}$ and ${\cal Y}
_{(h;\,k^{\prime}, \beta^{\prime})}$ in the function algebra 
$Fun_{h}\left(SL(2)\right)$, as obtained in ~(\ref{eq:xembed}) and
~(\ref{eq:yembed}) respectively, by using the infinitesimal 
invariance of the genarating elements with respect to the 
appropriate twisted primitive elements of the corresponding 
dual universal enveloping algebra $U_{h}(sl(2))$. This approach
was first developed by Koornwinder ~\cite{K93} in the context
of the Podles quantum sphere ~\cite{P87} for the standard
$q$-deformed function algebra $Fun_{q}\left(SL(2)\right)$. 
For our purpose, we first recall certain definitions
~\cite{K93} regarding dual Hopf algebras. 

Two Hopf algebras ${\cal U}$ and ${\cal A}$ are said 
to be in duality if there is a doubly nondegenerate 
bilinear form $(\sf {u}, \sf {a})\mapsto \langle \sf{u}, 
\sf{a}\rangle$ such that for $(\sf {u}, \sf{v})\in 
\cal{U}, (\sf {a},\sf {b})\in \cal{A}$ we have 
\begin{eqnarray}
&&\langle\bigtriangleup({\sf u}), {\sf a} \otimes {\sf b}\rangle
= \langle {\sf u}, {\sf ab} \rangle,\qquad 
\langle{\sf u} \otimes {\sf v}, \bigtriangleup(a)\rangle
= \langle {\sf uv}, {\sf a} \rangle,\nonumber\\
&&\langle 1_{\cal{U}}, {\sf a}\rangle 
= \varepsilon_{\cal{A}}({\sf a}),
\quad \langle {\sf u}, 1_{\cal{A}}\rangle 
= \varepsilon_{\cal{U}}({\sf u}),
\quad \langle S({\sf u}), {\sf a}\rangle =
\langle {\sf u}, S({\sf a}) \rangle.
\label{eq:dualHopf}
\end{eqnarray}
Instead of $\langle {\sf u}, {\sf a} \rangle$ we will also write
${\sf u}({\sf a})$ or ${\sf a}({\sf u})$.
An element $\sf {u}$ of the Hopf algebra ${\cal U}$ is called 
twisted primitive with respect to a group-like element $\sf {g}$
of ${\cal U}$ if
\begin{equation}
\bigtriangleup({\sf u}) = {\sf u}\otimes {\sf g} +
{\sf g^{-1}} \otimes {\sf u},\quad \varepsilon ({\sf u}) = 0,
\quad S({\sf u}) = -\,{\sf g u g^{-1}}.
\label{eq:twistprim}
\end{equation}
These twisted primitive elements play a key role in our 
present argument. For the elements ${\sf u}\in {\cal U}$ 
and ${\sf a}\in {\cal A}$, where ${\cal U}$ and ${\cal A}$ are
Hopf algebras in duality, we define elements ${\sf a.u}$ 
and ${\sf u.a}$ of ${\cal A}$ by
\begin{equation}
{\sf u.a} = ({\hbox {id}} \otimes {\sf u})\left(\bigtriangleup
({\sf a})\right),\qquad 
{\sf a.u} = ({\sf u} \otimes {\hbox {id}})\left(\bigtriangleup
({\sf a})\right).
\label{eq:Hopfdot}
\end{equation}
The operations defined in ~(\ref {eq:Hopfdot}) are left, 
respectively, right algebra actions of ${\cal U}$ of ${\cal A}$.
Assuming the standard notation 
\begin{equation}
\bigtriangleup ({\sf u}) =\sum_{(\sf {u})} {\sf u_{(1)}}
\otimes{\sf u_{(2)}},  
\label{Sweedler}
\end{equation}
we now obtain
\begin{equation}
{\sf u}.({\sf ab}) = \sum_{(\sf {u})} ({\sf u_{(1)}.a})
({\sf u_{(2)}.b}),\quad
({\sf ab}).{\sf u} = \sum_{(\sf {u})} ({\sf  a.u_{(1)}})
({\sf b.u_{(2)}})
\label{eq:algDot}
\end{equation}
for the elements ${\sf u}\in {\cal U}$ and $({\sf a,b})\in {\cal A}$.
We call an element ${\sf a}\in {\cal A}$ left (right) invariant
with respect to an element ${\sf u}\in {\cal U}$ if we have
${\sf u.a} = \varepsilon ({\sf u})\,{\sf a}$, respectively,
${\sf a.u} = \varepsilon ({\sf u})\,{\sf a}$. The unit element of 
${\cal A}$ is bi-invariant with respect to all 
${\sf u}\in {\cal U}$. If ${\sf u}$ is a twisted primitive 
element, then $\varepsilon ({\sf u}) = 0$, and it follows from
~(\ref{eq:algDot}) 
\begin{eqnarray}  
{\sf u.a} = 0 \quad {\hbox {and}} \quad {\sf u.b} = 0 \quad
\Longrightarrow \quad {\sf u.(ab)} = 0,\nonumber\\
{\sf a.u} = 0 \quad {\hbox {and}} \quad {\sf b.u} = 0 \quad
\Longrightarrow \quad {\sf (ab).u} = 0.
\label{eq:invarAlg}
\end{eqnarray}
As evidenced from (\ref{eq:invarAlg}) the left (right) invariant 
elements of ${\cal A}$ with respect to some twisted primitive 
element of ${\cal U}$ form a unital subalgebra of ${\cal A}$.   

In our context, the Jordanian  Hopf algebras $U_{h}(sl(2))$ and
$Fun_{h}\left(SL(2)\right)$ are in duality with the following
pairing of the generators:
\begin{eqnarray}
&&\left\langle T^{\pm 1},\left(\begin{array}{cc}a & b\\
c & d \end{array}\right)  \right\rangle = \left(\begin{array}{cc}
1 & {\pm} h\\0 & 1\end{array}\right),\quad
\left\langle Y,\left(\begin{array}{cc}a & b\\
c & d \end{array}\right)
\right\rangle = \left(\begin{array}{cc}
0 & 0\\1 & 0\end{array}\right),\nonumber\\
&&\left\langle H,\left(\begin{array}{cc}a & b\\
c & d \end{array}\right) \right\rangle = 
\left(\begin{array}{cc}1 & 0\\0 & -1\end{array}\right).
\label{eq:pairing}
\end{eqnarray}
{}From the coalgebraic structure ~(\ref{eq:jcoalg}) of the Hopf
algebra $U_{h}(sl(2))$ we see that any linear combination 
of $(T - T^{-1}), Y$ and $H$ is twisted primitive with respect
to the group-like element $T$ of the universal enveloping
algebra $U_{h}(sl(2))$.  
The embedding of the algebras ${\cal X}_{(h;\; k, \beta)}$ 
and ${\cal Y}_{(h;\; k^{\prime}, \beta^{\prime)}}$ in the function 
algebra $Fun_{h}\left(SL(2)\right)$ may now be naturally obtained
using the previously described invariance properties with respect
to these twisted primitive elements.
The generators $\left({\sf x_{i}\, | i= -1, 0, 1} \right)$ 
of the algebra ${\cal X}_{(h;\; k, \beta)}$ are left invariant 
with respect to a twisted primitive element of the dual 
algebra $U_{h}(sl(2))$:
\begin{equation}
{\cal P_{L}}\,.\,{\sf x_{j}} = 0 \quad {\hbox {for}}\quad
{\sf {j=(-1,0,1)}},
\label{eq:hlftinvr}
\end{equation}
where
\begin{equation}
{\cal P_{L}} = k {\rho}^{-1}\, \left(1 + \frac{3}{2}\,h^2\right)\,
\left(\frac{T-T^{-1}}{2h}\right) - H + 2\,k {\rho}^{-1}\, Y.
\label{eq:xprim}
\end{equation}  
Similarly the generators $\left({\sf y_{i}\, | i= -1, 0, 1} \right)$
of the algebra ${\cal Y}_{(h;\; k^{\prime}, \beta^{\prime)}}$
are right invariant with respsect to {\it another} twisted primitive 
element ${\cal P_{R}}$ of the $U_{h}(sl(2))$ algebra. This 
leads to an embedding of the  ${\cal Y}_{(h;\; k^{\prime}, \beta
^{\prime)}}$ algebra $Fun_{h}\left(SL(2)\right)$. The  above 
condition of right invarinnce reads   
\begin{equation}
{\sf y_{j}}\,.\, {\cal P_{R}} = 0 \quad {\hbox {for}}\quad
{\sf {j=(-1,0,1)}},
\label{eq:hrhtinvr}
\end{equation}
where the corresponding twisted primitive element is
\begin{equation}
{\cal P_{R}} = k^{\prime}{\rho^{\prime}}^{\,-1}\, 
\left(1 - \half\,h^2\right)\,
\left(\frac{T-T^{-1}}{2h}\right) - H 
+ 2\,k^{\prime} {\rho^{\prime}}^{\,-1}\, Y.
\label{eq:yprim}
\end{equation}
It is interesting to point out that the twisted primitive elements
${\cal P_{L}}$ and ${\cal P_{R}}$ implementing the left  
and the right invariances of the algebras ${\cal X}_{(h;\; k, \beta)}$
and ${\cal Y}_{(h;\; k^{\prime}, \beta^{\prime)}}$ respectively, 
are , in general,{\it distinct}. Except of the special values
of the parameters $k \rho^{-1} \rightarrow 0$,
$k^{\prime}{\rho^{\prime}}^{\,-1} \rightarrow 0$, the elements
${\cal P_{L}}$ and ${\cal P_{R}}$ defined by (\ref{eq:xprim}) and
(\ref{eq:yprim}) respectively, differ for a nonzero value of the
deformation parameter $h$. These twisted primitives, however,
become identical  in the $\rho \rightarrow \infty$
and $\rho^{\prime} \rightarrow \infty$ limits discussed earlier
in the context of equations (\ref{eq:xinfembed}) and 
(\ref{eq:yinfembed}) respectively. In these limits we obtain
\begin{equation}
\left({\cal P_{L}}\right)_{\rho \rightarrow \infty}\,\sim H,
\qquad \left({\cal P_{R}}\right)_{\rho^{\prime} \rightarrow 
\infty}\,\sim H.
\label{eq:infprim}
\end{equation}

\sect{Conclusion}
In conclusion, we have constructed a one parameter family of
algebras invariant under the left (right) coactions of the
group-like elements ${\cal T}_{h}^{(j=1)}$ of the Jordanian
function algebra $Fun_{h}\left(SL(2)\right)$. These algebras
may be regarded as Jordanian quantization of the two-dimensional 
spheres. Algebraically, an invertible map exist between the quantum 
spaces invariant under the left and the right coactions of the
${\cal T}_{h}^{(j=1)}$ matrix. The embedding of these algebras
in the function algebra $Fun_{h}\left(SL(2)\right)$ may be understood 
as an infinitesimal invariance of the generating elements of
these algebras with respect to the appropriate twisted primitive
elements of the Jordanian enveloping algebra $U_{h}(sl(2))$. It
is interesting to point out that , unlike the standard 
$q$-deformed case~\cite{K93}, in the present Jordanian example
the twisted primitive elements ${\cal P_{L}}$ and ${\cal P_{R}}$
associated with the spaces invariant under the left and the 
coactions of the ${\cal T}_{h}^{(j=1)}$ matrix respectively, are
distinct except for special values of the parameters.

As applications and further extensions of the present work, we
wish consider the followings. The classification of the differential
calculi on the Jordanian deformed sphere may be obtained. We will  
present the results elsewhere. In another development, the authors
of the Ref. \cite{MS00} considered the invariant propagator
of the $h$-deformed Laplacian on the $h$-deformed Lobachevski 
plane. Now this propagator may be studied on the Jordanian spherical
spaces. Lastly, we have recently constructed ~\cite{CC00}
a Jordanian deformed quasi-Hopf enveloping algebra satisfying
a shifted Yang-Baxter equation. The present problem may be 
investigated in the framework of Jordanian quasi-Hopf
function algebra.   

\noindent{\bf Acknowledgement} 

\bigskip 
 
\noindent We wish to thank Amitabha Chakrabarti for fruitful 
correspondence.    

\bibliographystyle{amsplain}

\end{document}